\documentclass{amsart}

\usepackage[margin=2cm]{geometry}
\usepackage{amsmath,amssymb,amsthm,mathtools}
\usepackage{enumerate}
\PassOptionsToPackage{hyphens}{url}
\usepackage{xcolor}
\usepackage[colorlinks=true,linkcolor=blue,citecolor=blue,urlcolor=blue]{hyperref}

\newtheorem{theorem}{Theorem}[section]
\newtheorem{lemma}[theorem]{Lemma}

\newtheorem{corollary}[theorem]{Corollary}
\theoremstyle{definition}

\newtheorem{remark}[theorem]{Remark}
\newtheorem{example}[theorem]{Example}

\newcommand{\N}{\mathbb N}
\newcommand{\Z}{\mathbb Z}

\newcommand{\R}{\mathbb R}
\newcommand{\PP}{\mathcal P}
\newcommand{\E}{\mathbb E}
\newcommand{\one}{\mathbf 1}
\newcommand{\abs}[1]{\left|#1\right|}
\newcommand{\norm}[1]{\left\|#1\right\|}
\newcommand{\floor}[1]{\left\lfloor #1\right\rfloor}
\newcommand{\muG}{\mu_G}

\numberwithin{equation}{section}

\begin{document}

\title[Khinchin constants along deterministic subsequences]{Khinchin's constant along deterministic subsequences} 

\author[M. Hussain]{Mumtaz Hussain}
\address{Mumtaz Hussain, Department of Mathematical and Physical Sciences, La Trobe University, Bendigo 3552, Australia}
\email{m.hussain@latrobe.edu.au}

\date{}

\begin{abstract}
Khinchin's theorem says that, for Lebesgue almost every real number, the geometric mean of the partial quotients in its continued fraction expansion converges to a universal constant $K\approx 2.6854520010\ldots$.  We show that the same constant is obtained when the partial quotients are sampled along any fixed strictly increasing deterministic sequence of indices.  Thus arithmetic progressions, polynomial sequences, prime indices, and lacunary sequences all give the usual Khinchin constant.  We also look at a different selection rule: instead of restricting the indices, we restrict the values of the partial quotients.  In this case the answer changes.  The geometric mean of the prime-valued partial quotients converges almost surely to a new constant $K_{\PP}\approx {3.439477971}\ldots$.
\end{abstract}


\maketitle

\section{Introduction}
Let $        x=[0;a_1(x),a_2(x),a_3(x),\ldots]$ 
be the continued fraction expansion of an irrational number $x\in[0,1)$.  For a fixed irrational $x$, the partial quotients may behave very irregularly.  Nevertheless, {Khinchin's theorem \cite{Khinchin}} gives the following universal law: for Lebesgue almost every $x$,
\begin{equation*}\label{eq:classical-Khinchin}
        \lim_{N\to\infty}
        \left(\prod_{n=1}^{N}a_n(x)\right)^{1/N}
        =K,
\end{equation*}
where
\begin{equation}\label{eq:K-def}
{
        K
        =
        \prod_{r=1}^{\infty}
        \left(1+\frac1{r(r+2)}\right)^{\frac{\log r}{\log 2}}
        \approx 2.6854520010\ldots.
}
\end{equation}
The usual proof applies Birkhoff's ergodic theorem to the Gauss map
{
\[
        T(0)=0,
        \qquad
        T(x)=\frac1x-\left\lfloor\frac1x\right\rfloor
        \quad (x\in(0,1)).
\]
}
and to the integrable function $x\mapsto \log a_1(x)$.  Indeed,
\[
        a_n(x)=a_1(T^{n-1}x),
\]
and $T$ preserves the Gauss measure
\[
        d\muG(x)=\frac{1}{(1+x)\log 2}\,dx,
\]
which is equivalent to Lebesgue measure.

The purpose of this note is to distinguish between two natural ways of selecting partial quotients.  The first is to select their positions.  Fix a strictly increasing sequence of positive integers
\[
        1\le t_1<t_2<t_3<\cdots
\]
and consider
\[
        \left(\prod_{n=1}^{N}a_{t_n}(x)\right)^{1/N}.
\]
Our first theorem says that, for each fixed choice of $(t_n)$, the limiting constant is still $K$.  In other words, deterministic sampling of the indices does not create a new Khinchin constant.

\noindent This question is partly motivated by recent work on arithmetic progressions in the Gauss system.  Products such~as
\[
        a_n(x)a_{2n}(x)\cdots a_{dn}(x)
\]
appear naturally in connection with exceptional sets for the convergence criterion of multiple ergodic averages of Bourgain (1988) and Host--Kra (2005).  For the Gauss dynamical system, failure of these convergence criteria can be expressed through unusually large products of partial quotients along arithmetic progressions; see for instance \cite{Bourgain89,HostKra,HLSAMH, HussainShulgaAP}.  The result below shows that, as long as the indexing sequence is fixed in advance, the ordinary metric constant survives: arithmetic progressions, polynomial sequences, prime indices, and lacunary sequences all lead back to Khinchin's constant. {Here ``arithmetic progression'' means a fixed subsequence such as $t_n=mn+b$.  The triangular arrays $n,2n,\ldots,dn$, in which the sampled indices themselves depend on the row parameter, are not covered directly by the theorem; see Remark~\ref{rem:ap-arrays}.} 

Questions concerning ergodic averages along deterministic subsequences in the Gauss dynamical system have a long history. For instance, early contributions of Nair
\cite{Nair1991,Nair1994} studied continued fractions together with
ergodic averages along polynomial and other deterministic subsequences, 
establishing refinements of Birkhoff's ergodic theorem for deterministic sampling problems.  Although these works do not consider Khinchin's constant itself, they provide part of the historical motivation for studying deterministic sampling of partial quotients.

The second selection rule is different.  Instead of selecting positions, we select digit values.  Let $\PP$ be the set of prime numbers.  We keep only those partial quotients for which $a_n(x)\in\PP$. {Prime digits in typical continued fraction expansions have recently been studied in detail by Schindler and Zweim\"uller~\cite{SchindlerZweimueller}.  The prime-valued geometric mean considered here is a natural conditional ergodic statistic, included to contrast value selection with index selection.} 
The geometric mean over these digits no longer converges to $K$; it converges to the conditional constant obtained from the Gauss digit distribution after conditioning on the event that the digit is prime.  More precisely, if
\[
        N_{\PP}(x,N)=\#\{1\le n\le N:a_n(x)\in\PP\},
\]
then, for almost every $x$, {the following expression is defined for all sufficiently large $N$, and}
\[
        \left(\prod_{\substack{1\le n\le N\\ a_n(x)\in\PP}}a_n(x)\right)^{1/N_{\PP}(x,N)}
        \longrightarrow K_{\PP}       =
        \exp\left(
        \frac{
        \displaystyle\sum_{p\in\PP}
        \log p\,\log\left(1+\frac1{p(p+2)}\right)
        }{
        \displaystyle\sum_{p\in\PP}
        \log\left(1+\frac1{p(p+2)}\right)
        }
        \right)
        \approx{3.439477971}\ldots.
\]
Thus prime indices and prime-valued partial quotients should not be confused: prime indices give the classical constant $K$, while prime-valued digits give a genuinely different constant  $K_{\PP}$.

\medskip
\noindent
{The key technical ingredient in the proof of Theorem~\ref{thm:deterministic} is a strong law for sampled digit functions (Lemma~\ref{lem:slln-subseq} below). This lemma, which may be of independent interest, shows that any sequence of digits sampled along a fixed increasing sequence of indices obeys the usual law of large numbers, provided the function has finite fourth moment and the digit process is exponentially mixing.}

We now formally state the results.

\begin{theorem}\label{thm:deterministic}
Let $(t_n)_{n\ge1}$ be a fixed strictly increasing sequence of positive integers.  Then, for Lebesgue almost every $x\in[0,1)$,
\begin{equation*}\label{eq:deterministic-limit}
        \lim_{N\to\infty}
        \left(\prod_{n=1}^{N}a_{t_n}(x)\right)^{1/N}
        =K,
\end{equation*}
where $K$ is Khinchin's constant defined in \eqref{eq:K-def}.  Equivalently,
\[
        \lim_{N\to\infty}
        \frac1N\sum_{n=1}^{N}\log a_{t_n}(x)
        =\log K.
\]
\end{theorem}

\begin{example}\label{ex:examples}
Theorem \ref{thm:deterministic} applies to every fixed strictly increasing deterministic sequence of indices.  In particular, for Lebesgue almost every $x$:
\begin{enumerate}[(i)]
\item if $t_n=mn+b$ with $m\ge1$ and $b\ge0$, then
\[
        \left(\prod_{n=1}^{N}a_{mn+b}(x)\right)^{1/N}\longrightarrow K;
\]
\item if $P\in\Z[t]$ is eventually positive and strictly increasing, then, after discarding finitely many initial terms if necessary,
\[
        \left(\prod_{n=1}^{N}a_{P(n)}(x)\right)^{1/N}\longrightarrow K;
\]
\item if $p_n$ denotes the $n$-th prime number, then
\[
        \left(\prod_{n=1}^{N}a_{p_n}(x)\right)^{1/N}\longrightarrow K;
\]
\item if $t_n=\floor{\lambda^n}$ with $\lambda>1$, then the same conclusion holds after deleting the finitely many possible repetitions at the beginning.  More generally, any fixed strictly increasing lacunary sequence gives
\[
        \left(\prod_{n=1}^{N}a_{t_n}(x)\right)^{1/N}\longrightarrow K.
\]
\end{enumerate}
\end{example}

\begin{theorem}\label{thm:prime-valued}
Let $\PP$ be the set of prime numbers and define
\[
        N_{\PP}(x,N)=\#\{1\le n\le N:a_n(x)\in\PP\}.
\]
Then, $N_{\PP}(x,N)\to\infty$ for Lebesgue almost every $x$, and
\begin{equation}\label{eq:prime-limit}
        \lim_{N\to\infty}
        \left(
        \prod_{\substack{1\le n\le N\\ a_n(x)\in\PP}}
        a_n(x)
        \right)^{1/N_{\PP}(x,N)}
        =K_{\PP},
\end{equation}
where
\begin{equation}\label{eq:KP-def}
        K_{\PP}
        =
        \exp\left(
        \frac{
        \displaystyle\sum_{p\in\PP}
        \log p\,\log\left(1+\frac1{p(p+2)}\right)
        }{
        \displaystyle\sum_{p\in\PP}
        \log\left(1+\frac1{p(p+2)}\right)
        }
        \right)
    \approx {3.439477971}\ldots .
\end{equation}
\end{theorem}

\begin{corollary}\label{cor:prime-subseq}
Let $(t_n)_{n\ge1}$ be a fixed strictly increasing sequence of positive integers and put
\[
        N_{\PP,t}(x,N)=\#\{1\le n\le N:a_{t_n}(x)\in\PP\}.
\]
Then $N_{\PP,t}(x,N)\to\infty$ for Lebesgue almost every $x$, and
\[
        \left(
        \prod_{\substack{1\le n\le N\\ a_{t_n}(x)\in\PP}}
        a_{t_n}(x)
        \right)^{1/N_{\PP,t}(x,N)}
        \longrightarrow K_{\PP}.
\]
\end{corollary}

{
\begin{remark}\label{rem:exceptional-set}
The sequence $(t_n)$ is fixed before $x$ is chosen.  The exceptional null set in Theorem~\ref{thm:deterministic} may depend on $(t_n)$; the theorem is not a simultaneous assertion over all strictly increasing sequences.
\end{remark}
}

\section{Gauss measure and a strong law for sampled digits}
{The Gauss measure $\mu_G$ is invariant and ergodic for $T$, and it is equivalent to Lebesgue measure.  Therefore it is enough to prove all almost-everywhere statements with respect to $\muG$.}

For $r\in\N$, the one-cylinder $\{a_1=r\}$ is the interval $(1/(r+1),1/r]$.  Hence
\begin{equation*}\label{eq:digit-distribution}
        \muG(a_1=r)
        =\frac1{\log 2}\int_{1/(r+1)}^{1/r}\frac{dx}{1+x}
        =\frac1{\log 2}\log\frac{(r+1)^2}{r(r+2)}
        =\frac1{\log 2}\log\left(1+\frac1{r(r+2)}\right).
\end{equation*}
Consequently,
\begin{equation}\label{eq:int-log-a}
        \int\log a_1\,d\muG
        =
        \sum_{r=1}^{\infty}
        \frac{\log r}{\log 2}
        \log\left(1+\frac1{r(r+2)}\right)
        =\log K.
\end{equation}
The integral is finite, since the summand is comparable to $(\log r)/r^2$.

Similarly,
\begin{equation}\label{eq:prime-prob}
        \muG(a_1\in\PP)
        =
        \sum_{p\in\PP}
        \frac1{\log 2}
        \log\left(1+\frac1{p(p+2)}\right)
        >0,
\end{equation}
and
\begin{equation}\label{eq:int-prime-log}
        \int \log a_1\,\one_{\{a_1\in\PP\}}\,d\muG
        =
        \sum_{p\in\PP}
        \frac{\log p}{\log 2}
        \log\left(1+\frac1{p(p+2)}\right).
\end{equation}
Both series converge because
\[
        \log\left(1+\frac1{p(p+2)}\right)\asymp p^{-2}.
\]


\medskip

{
For integers $1\le u\le v\le\infty$, let
\[
        \mathcal F_u^v=\sigma(a_u,a_{u+1},\ldots,a_v),
\]
with the evident interpretation when $v=\infty$.  Under the Gauss measure, the digit process $(a_n)_{n\ge1}$ is strictly stationary and exponentially $\psi$-mixing: there are constants $C_0>0$ and $0<\rho<1$ such that, for every $i,k\ge1$, every $A\in\mathcal F_1^i$, and every $B\in\mathcal F_{i+k}^{\infty}$,
\begin{equation}\label{eq:psi-mixing}
        \abs{\muG(A\cap B)-\muG(A)\muG(B)}
        \le C_0\rho^k\muG(A)\muG(B).
\end{equation}
This is classical; see, for example, \cite[Section~2.3.4]{IosifescuKraaikamp}, \cite[Lemma~2.2]{KessebohmerSchindler}, or \cite[Section~4]{SchindlerZweimueller}.  By stationarity, the same estimate holds for two arbitrary digit blocks separated by at least $k$ positions.  Throughout this section, $\E$ denotes expectation with respect to $\muG$.

We shall repeatedly use the following consequence of \eqref{eq:psi-mixing}.  If $U$ and $V$ are integrable random variables measurable with respect to two digit blocks separated by at least $k$ positions, then
\begin{equation}\label{eq:functional-psi}
        \abs{\E(UV)-\E(U)\E(V)}
        \le C_0\rho^k\E\abs{U}\,\E\abs{V}.
\end{equation}
 Indeed, for non-negative simple functions this follows by expanding them as finite linear combinations of indicator functions and applying \eqref{eq:psi-mixing} term by term.  Monotone convergence gives the estimate for non-negative integrable functions, and decomposing into positive and negative parts gives the general case.

\begin{lemma}\label{lem:fourth-moment}
Let $F:\N\to\R$ satisfy
\[
        \sum_{r=1}^{\infty} |F(r)|^4\,\muG(a_1=r)<\infty.
\]
Put
\[
        \varphi_F(x)=F(a_1(x)),
        \qquad
        \overline\varphi_F=\varphi_F-\int\varphi_F\,d\muG.
\]
Then there is a constant $C_F>0$, independent of the sequence $(t_n)$, such that for every strictly increasing sequence $(t_n)$ of positive integers and every $N\ge1$,
\begin{equation}\label{eq:fourth-moment}
        \int
        \left|
        \sum_{n=1}^{N}
        \overline\varphi_F(T^{t_n-1}x)
        \right|^4
        d\muG(x)
        \le C_F N^2.
\end{equation}
\end{lemma}

\begin{proof}
Define the centered stationary process
\[
        X_m(x)=\overline\varphi_F(T^{m-1}x),\qquad m\ge1,
\]
and put $Y_n=X_{t_n}$.  Set $L=\norm{X_1}_{L^4(\muG)}<\infty$.  Expanding the fourth power and then taking absolute values gives
\[
        \E\left|\sum_{i=1}^{N}Y_i\right|^4
        \le
        \sum_{i,j,k,\ell=1}^{N}
        \abs{\E(Y_iY_jY_kY_\ell)}.
\]
After ordering the four indices, each non-decreasing quadruple
\[
        1\le i_1\le i_2\le i_3\le i_4\le N
\]
occurs at most $4!$ times.  Write
\[
        h_1=i_2-i_1,\qquad h_2=i_3-i_2,\qquad h_3=i_4-i_3,
        \qquad h=\max\{h_1,h_2,h_3\},
\]
and, when the maximum is attained more than once, choose its leftmost occurrence.  Fix an integer $M\ge1$.

If $h\le M$, then, for each $i_1$, there are $O_M(1)$ choices for the remaining indices.  Hence there are $O_M(N)$ such quadruples, and H\"older's inequality gives
\[
        \abs{\E(Y_{i_1}Y_{i_2}Y_{i_3}Y_{i_4})}\le L^4.
\]
Their total contribution is therefore $O_M(N)$.

Suppose now that $h>M$.  If the chosen largest gap is $h_1$, then $Y_{i_1}$ and $Y_{i_2}Y_{i_3}Y_{i_4}$ are measurable with respect to digit blocks separated by at least
\[
        t_{i_2}-t_{i_1}\ge i_2-i_1=h.
\]
Since $\E Y_{i_1}=0$, \eqref{eq:functional-psi} and H\"older's inequality yield
\[
\begin{split}
        \abs{\E(Y_{i_1}Y_{i_2}Y_{i_3}Y_{i_4})}
        &\le C_0\rho^h\E\abs{Y_{i_1}}\,
        \E\abs{Y_{i_2}Y_{i_3}Y_{i_4}}\\
        &\le C_0\rho^h L^4.
\end{split}
\]
The case in which the chosen largest gap is $h_3$ is identical.

If the chosen largest gap is $h_2$, put
\[
        A=Y_{i_1}Y_{i_2},\qquad B=Y_{i_3}Y_{i_4}.
\]
Then \eqref{eq:functional-psi} gives
\begin{equation}\label{eq:middle-gap}
        \abs{\E(AB)}
        \le
        \abs{\E(A)}\,\abs{\E(B)}+C_0\rho^h L^4,
\end{equation}
where we used $\E|A|,\E|B|\le L^2$.

For any fixed $h$, and for each of the three possible positions of the chosen largest gap, the number of non-decreasing quadruples is $O(N(h+1)^2)$: one chooses the location of the gap and the two remaining gaps, each of which is at most $h$.  Consequently, all the exponentially decaying terms above contribute at most
\[
        O\!\left(N\sum_{h>M}(h+1)^2\rho^h\right)=O(N).
\]

It remains to sum the product term in \eqref{eq:middle-gap}.  Dropping the ordering and largest-gap restrictions only enlarges the sum, so
\[
\begin{split}
&\sum_{1\le i_1\le i_2\le i_3\le i_4\le N}
        \abs{\E(Y_{i_1}Y_{i_2})}\,
        \abs{\E(Y_{i_3}Y_{i_4})}\le
        \left(
        \sum_{1\le i\le j\le N}
        \abs{\E(Y_iY_j)}
        \right)^2.
\end{split}
\]
For $i<j$, the variables are centered and their digit positions are separated by at least $t_j-t_i\ge j-i$.  Hence \eqref{eq:functional-psi} gives
\[
        \abs{\E(Y_iY_j)}
        \le C_0\rho^{t_j-t_i}\E|Y_i|\E|Y_j|
        \le C_0L^2\rho^{j-i},
\]
while $\E(Y_i^2)\le L^2$.  It follows that
\[
        \sum_{1\le i\le j\le N}\abs{\E(Y_iY_j)}
        \le L^2N+C_0L^2\sum_{r=1}^{N-1}(N-r)\rho^r
        \ll_F N.
\]
Thus the product terms contribute $O_F(N^2)$.  Combining the preceding estimates proves \eqref{eq:fourth-moment}.
\end{proof}
}

\begin{lemma}\label{lem:slln-subseq}
Let $F:\N\to\R$ satisfy
\[
        \sum_{r=1}^{\infty}|F(r)|^4\muG(a_1=r)<\infty.
\]
Then, for every fixed strictly increasing sequence $(t_n)_{n\ge1}$ of positive integers,
\begin{equation}\label{eq:slln-subseq}
        \frac1N\sum_{n=1}^{N}
        \left(F(a_{t_n}(x))-\int F(a_1)\,d\muG\right)
        \longrightarrow0
\end{equation}
for $\muG$-almost every $x$.
\end{lemma}

\begin{proof}
Let
\[
        Z_n(x)=F(a_{t_n}(x))-\int F(a_1)\,d\muG.
\]
By Lemma \ref{lem:fourth-moment},
\[
        \E\left|\sum_{n=1}^{N}Z_n\right|^4\le C_F N^2.
\]
Therefore, for every $\varepsilon>0$, Markov's inequality gives
\[
        \muG\left(
        \left|\frac1N\sum_{n=1}^{N}Z_n\right|>\varepsilon
        \right)
        \le
        \frac{C_F}{\varepsilon^4 N^2}.
\]
The right-hand side is summable in $N$.  The Borel--Cantelli lemma implies that, for $\muG$-almost every $x$,
\[
        \left|\frac1N\sum_{n=1}^{N}Z_n(x)\right|>\varepsilon
\]
only for finitely many $N$.  Taking $\varepsilon$ over a countable sequence tending to zero proves \eqref{eq:slln-subseq}.
\end{proof}

\section{Proof of Theorem \ref{thm:deterministic}}
\begin{proof}
Apply Lemma \ref{lem:slln-subseq} with $
        F(r)=\log r.$ 
This choice is admissible because
\[
        \sum_{r=1}^{\infty}(\log r)^4\muG(a_1=r)<\infty,
\]
as $\muG(a_1=r)\asymp r^{-2}$.  Hence, for $\muG$-almost every $x$,
\[
        \frac1N\sum_{n=1}^{N}
        \left(\log a_{t_n}(x)-\int\log a_1\,d\muG\right)
        \longrightarrow0.
\]
By \eqref{eq:int-log-a}, we have
\[
        \int\log a_1\,d\muG=\log K.
\]
Therefore
\[
        \frac1N\sum_{n=1}^{N}\log a_{t_n}(x)
        \longrightarrow \log K
\]
for $\muG$-almost every $x$.  Since $\muG$ and Lebesgue measure are equivalent, the same assertion holds for Lebesgue almost every $x$.  Exponentiating gives
\[
        \left(\prod_{n=1}^{N}a_{t_n}(x)\right)^{1/N}
        \longrightarrow K.
\]
\end{proof}

\begin{remark}\label{rem:no-new-constant}
The expectation of each sampled logarithmic digit is always the same:
\[
        \int\log a_{t_n}\,d\muG
        =
        \int\log a_1\,d\muG
        =\log K.
\]
The only real point is to show that the fluctuations average out along the chosen deterministic sequence.  Exponential mixing gives exactly this.  The fact that $(t_n)$ is increasing only helps, because the sampled digits are at least as far apart as ordinary consecutive digits in the sampled order.
\end{remark}

\begin{remark}\label{rem:ap-arrays}
Theorem \ref{thm:deterministic} treats fixed subsequence averages, for instance $t_n=an+b$.  The row-dependent averages
\[
        \frac1n\sum_{i=1}^{n}\log a_{in}(x),
\]
which arise when studying products such as $        a_n(x)a_{2n}(x)\cdots a_{n^2}(x),$ 
are triangular rather than fixed-subsequence averages.  They are not covered directly by Theorem \ref{thm:deterministic}.  However, the same fourth-moment method is well suited to them, since within the $n$-th row the gaps between the indices $n,2n,\ldots,n^2$ are at least $n$, so the variables in one row become increasingly decorrelated.  The theorem above isolates the simpler fixed-subsequence phenomenon: once the indices are fixed in advance and strictly increasing, the limiting constant is the usual Khinchin constant.
\end{remark}

\section{Prime-valued partial quotients}
{We now prove the result for prime-valued digits.  The asymptotic frequency of prime digits is already recorded in \cite[Proposition~2.1]{SchindlerZweimueller}; the argument below additionally identifies the conditional geometric mean and extends it to every fixed deterministic subsequence.}  Define
\[
        f(x)=\one_{\{a_1(x)\in\PP\}},
        \qquad
        g(x)=\log a_1(x)\,\one_{\{a_1(x)\in\PP\}}.
\]
Both functions have finite fourth moment with respect to $\muG$.  This is clear for $f$, and for $g$ it follows from
\[
        \int |g|^4\,d\muG
        =
        \sum_{p\in\PP}(\log p)^4\muG(a_1=p)
        <\infty.
\]

\begin{proof}[Proof of Theorem \ref{thm:prime-valued}]
Apply Lemma \ref{lem:slln-subseq} to the full sequence $t_n=n$ and to $
        F(r)=\one_{\{r\in\PP\}}.$ 
We get
\begin{equation}\label{eq:prime-count-limit}
        \frac{N_{\PP}(x,N)}{N}
        =
        \frac1N\sum_{n=1}^{N}\one_{\{a_n(x)\in\PP\}}
        \longrightarrow
        \muG(a_1\in\PP)
\end{equation}
for $\muG$-almost every $x$.  By \eqref{eq:prime-prob}, this limit is positive.  Hence $N_{\PP}(x,N)\to\infty$ and, in particular, the denominator in \eqref{eq:prime-limit} is eventually nonzero for almost every $x$.

Next apply Lemma \ref{lem:slln-subseq} with $        F(r)=\log r\,\one_{\{r\in\PP\}}.$ 
This gives
\begin{equation}\label{eq:prime-log-sum-limit}
        \frac1N
        \sum_{\substack{1\le n\le N\\a_n(x)\in\PP}}
        \log a_n(x)
        \longrightarrow
        \int \log a_1\,\one_{\{a_1\in\PP\}}\,d\muG.
\end{equation}
Dividing \eqref{eq:prime-log-sum-limit} by \eqref{eq:prime-count-limit} yields
\[
        \frac{
        \displaystyle\sum_{\substack{1\le n\le N\\a_n(x)\in\PP}}
        \log a_n(x)
        }{
        N_{\PP}(x,N)
        }
        \longrightarrow
        \frac{
        \displaystyle\int \log a_1\,\one_{\{a_1\in\PP\}}\,d\muG
        }{
        \displaystyle\muG(a_1\in\PP)
        }.
\]
Using \eqref{eq:prime-prob} and \eqref{eq:int-prime-log}, the limit is
\[
        \frac{
        \displaystyle\sum_{p\in\PP}
        \log p\,\log\left(1+\frac1{p(p+2)}\right)
        }{
        \displaystyle\sum_{p\in\PP}
        \log\left(1+\frac1{p(p+2)}\right)
        }.
\]
Exponentiating proves \eqref{eq:prime-limit} with $K_{\PP}$ as in \eqref{eq:KP-def}.  Finally, since $\muG$ is equivalent to Lebesgue measure, the same conclusion holds for Lebesgue almost every $x$.
\end{proof}

\begin{remark}\label{rem:product-form}
The constant in Theorem \ref{thm:prime-valued} can also be written as
\[
        K_{\PP}
        =
        \prod_{p\in\PP}
        p^{\muG(a_1=p)/\muG(a_1\in\PP)}.
\]
Thus $K_{\PP}$ is simply the geometric mean of the Gauss digit distribution conditioned on the digit being prime.  This explains why it differs from $K$: the classical Khinchin constant averages over all possible digits, while $K_{\PP}$ averages only over prime digits.
\end{remark}

\begin{proof}[Proof of Corollary \ref{cor:prime-subseq}]
Apply Lemma \ref{lem:slln-subseq} twice along the given sequence $(t_n)$: first with $
        F(r)=\one_{\{r\in\PP\}},$ 
and then with $        F(r)=\log r\,\one_{\{r\in\PP\}}.$ 
This gives, for $\muG$-almost every $x$,
\[
        \frac{N_{\PP,t}(x,N)}{N}\to \muG(a_1\in\PP)>0
\]
and
\[
        \frac1N
        \sum_{\substack{1\le n\le N\\a_{t_n}(x)\in\PP}}
        \log a_{t_n}(x)
        \to
        \int \log a_1\,\one_{\{a_1\in\PP\}}\,d\muG.
\]
Taking the ratio gives the same conditional mean as in the proof of Theorem \ref{thm:prime-valued}, and exponentiating gives $K_{\PP}$.  The equivalence of $\muG$ and Lebesgue measure then gives the stated Lebesgue almost-everywhere conclusion.
\end{proof}

\noindent{\bf Acknowledgements.} 
The author thanks Radhakrishnan (Kit) Nair, Gerardo  Gonz\'alez Robert,  Bixuan Li, and Nikita Shulga for useful discussions and comments on an earlier draft.

\noindent \textit{AI Declaration:}  The principal mathematical ideas, results, arguments, and proofs presented in this paper originate with the author while working on \cite{HLSAMH, HussainShulgaAP} and have been discussed with various colleagues longtime ago. ChatGPT5.6 Sol have been used to improve language and assisting the author in checking arguments that he has developed and written himself. In any case, the author take full responsibility for its correctness.





\begin{thebibliography}{99}

\bibitem{Bourgain89}
J.~Bourgain.
\newblock Pointwise ergodic theorems for arithmetic sets.
\newblock \emph{Inst. Hautes \'Etudes Sci. Publ. Math.}, 69:5--45, 1989.



\bibitem{HostKra}
B.~Host and B.~Kra.
\newblock Nonconventional ergodic averages and nilmanifolds.
\newblock \emph{Ann. of Math. (2)}, 161(1):397--488, 2005.

{
\bibitem{HLSAMH}
M.~Hussain, B.~Li, and N.~Shulga.
\newblock Averages of Gauss map iterations in arithmetic progressions.
\newblock \emph{Acta Math. Hungar.}, 178(2):392-424,2026.}

{
\bibitem{HussainShulgaAP}
M.~Hussain and N.~Shulga.
\newblock Metrical properties of finite product of partial quotients in arithmetic progressions.
\newblock \emph{Ann. Sc. Norm. Super. Pisa Cl. Sci. (5)}, 27(1):337--376, 2026.
}


\bibitem{IosifescuKraaikamp}
M.~Iosifescu and C.~Kraaikamp.
\newblock \emph{Metrical Theory of Continued Fractions}.
\newblock Kluwer Academic Publishers, Dordrecht, 2002.




{
\bibitem{KessebohmerSchindler}
M.~Kesseb\"ohmer and T.~I. Schindler.
\newblock Limit theorems for counting large continued fraction digits.
\newblock \emph{Lithuanian Math. J.}, 60(2):189--207, 2020.


\bibitem{Nair1991}
R.~Nair.
\newblock On polynomials in primes and J. Bourgain's circle method approach to ergodic theorems.
\newblock \emph{Ergodic Theory Dynam. Systems}, no. 3, 485–499, 1991.


\bibitem{Nair1994}
R.~Nair.
\newblock On the metrical theory of continued fractions.
\newblock \emph{Proc. Amer. Math. Soc.}, 120(4), 1041-1046, 1994.



\bibitem{SchindlerZweimueller}
T.~I. Schindler and R.~Zweim\"uller.
\newblock Prime numbers in typical continued fraction expansions.
\newblock \emph{Boll. Unione Mat. Ital.}, 16(2):259--274, 2023.
}

\bibitem{Khinchin}
A.~Ya. Khinchin.
\newblock \emph{Continued Fractions}.
\newblock Dover Publications, Mineola, NY, 1997.

\end{thebibliography}
\end{document}